\documentclass[11pt]{article}
\usepackage{amsmath,amssymb,fullpage}
\def\8{\infty}

\def\la{\langle}

\def\ra{\rangle}

\newtheorem {thm}{Theorem}[section]
\newtheorem {lem}[thm]{Lemma}
\newtheorem {cor}[thm]{Corollary}
\newtheorem {prop}[thm]{Proposition}
\newtheorem {df}[thm]{Definition}
\makeatletter

\@addtoreset{equation}{section}
\makeatother
\begin{document}
\title{\bf\boldmath
Derivation of a 
$BC_n$ elliptic summation formula \\
via the fundamental invariants
%\\{\large (preliminary version)}
}
\author{
{\sc Masahiko Ito}\footnote{
School of Science and Technology for Future Life,
Tokyo Denki University,
Tokyo 120-8551, Japan 
}
\ \ 
and {\sc Masatoshi Noumi}\footnote{
Department of Mathematics, Kobe University, 
Rokko, Kobe 657 8501, Japan
}
}
\date{}
\maketitle
\begin{abstract}
We give an alternative proof of an elliptic summation formula of type $BC_n$ 
by applying the fundamental $BC_n$ invariants to the study of 
Jackson integrals associated with the summation formula.
\end{abstract}



%%%%%%%%%%%%%%%%%%%%%%%%%%%%%%%%%%%%
%%%%%%%%%%%%%%%%%%%%%%%%%%%%%%%%%%%%
\section{Introduction}
\label{intro}
In this paper we investigate the following 
multiple elliptic summation formula of type $BC_n$:  
Under the conditions  
$a_1\cdots a_6\,t^{2n-2}=q$ and $a_1a_6\,t^{n-1}=q^{-N}$,   
\begin{eqnarray}\label{BCn-sum}
&&
%\hspace{-5pt}
\sum_{N\ge\nu_1\ge\cdots
\ge\nu_n\ge 0}%\!\!
\left[
\prod_{i=1}^{n}
(qt^{2(i-1)})^{\nu_i}
\frac{\theta(q^{2\nu_i}\xi_i^{2};p)}{\theta(\xi_i^{2};p)}
%\hspace{-2pt}
\prod_{1\le j<k\le n}
%\hspace{-4pt}
\frac
{\theta(q^{\nu_j-\nu_k}\xi_j/\xi_k;p)\theta(q^{\nu_j+\nu_k}\xi_j\xi_k;p)}
{\theta(\xi_j/\xi_k;p)\theta(\xi_j\xi_k;p)}
\right.
\nonumber\\[-2pt]
&&
\left.
\hspace{20pt}
\times
\prod_{i=1}^{n}
\prod_{m=1}^{6}
%\!
\frac{\theta(a_m\xi_i;p;q)_{\nu_i}}
{\theta(qa_m^{-1}\xi_i;p;q)_{\nu_i}}
%\!
\prod_{1\le j<k\le n}
%\!
\frac
{\theta(t\xi_j/\xi_k;p;q)_{\nu_j-\nu_k}
\theta(t\xi_j\xi_k;p;q)_{\nu_j+\nu_k}}
{\theta(qt^{-1}\xi_j/\xi_k;p;q)_{\nu_j-\nu_k}
\theta(qt^{-1}\xi_j\xi_k;p;q)_{\nu_j+\nu_k}\!}
\right]
\nonumber\\[5pt]
&&\hspace{10pt}=
\prod_{i=1}^{n}
\frac{
\theta(qa_1^2t^{n+i-2};p;q)_N
\prod_{2\le j<k\le 4}\theta(qa_j^{-1}a_k^{-1}t^{1-i};p;q)_N
}{
\theta(qa_1^{-1}a_2^{-1}a_3^{-1}a_4^{-1}t^{2-n-i};p;q)_N
\prod_{m=2}^{4}\theta(qa_m^{-1}a_1t^{i-1};p;q)_N
}
\end{eqnarray}
for $\xi_i=a_1t^{n-i}$ ($i=1,\ldots,n$), where 
we have used the standard notation 
%\pagebreak
\begin{equation*}
(x;p)_\infty=\prod_{i=0}^{\infty}(1-p^ix),\quad
\theta(x;p)=(x;p)_\infty(p/x;p)_\infty
\end{equation*}
for the shifted factorial and the theta function of base $p$ ($|p|<1$).  
We have also used the notation 
\begin{equation*}
\theta(x;p;q)_k=\theta(x;p)\theta(qx;p)\cdots\theta(q^{k-1}x;p)
\quad(k=0,1,2,\ldots),
\end{equation*}
for the elliptic shifted factorial as in \cite{vDS2000}. 
(We decided not to use the notation $(x;q,p)_k$ in \cite{GR2004,W2002} 
to avoid confusion with the infinite product $(x;p,q)_\infty$.)
\par
The summation formula \eqref{BCn-sum} for $n=1$ is the Frenkel--Turaev sum \cite{FT1997}.
This formula for general $n$ was first conjectured by Warnaar \cite{W2002}, and 
proved by van Diejen--Spiridonov \cite{vDS2000} in the case where $p=0$. 
The proof in the general case was given by Rosengren \cite{Ro2001}, 
Rains \cite{Ra2006} and Coskun--Gustafson \cite{CG2007}. 
(See also Schlosser \cite{Sc2007} in the special case $t=q$.) 
\par
The purpose of this paper is to give an alternative proof 
of the summation formula \eqref{BCn-sum} based 
on the idea of the Jackson integrals and the associated 
$q$-difference de Rham cohomology.   
In our approach 
a certain family of fundamental $BC_n$ invariants plays an important role 
in analyzing the structure of this summation. 
\par
After formulating the summation formula \eqref{BCn-sum} in terms of 
Jackson integrals in Section \ref{Sect.2}, we introduce in Section \ref{Sect.3} a family of holomorphic functions 
$E_r(z)$ ($r=0,1,\ldots,n$), 
called {\em fundamental invariants}, 
which are characterized by a certain interpolation property.  
In Section \ref{Sect.4} we establish two-term relations for the Jackson integrals involving $E_r(z)$.  
On the basis of these preparations, we give in Section \ref{Sect.5} a proof of the summation formula 
\eqref{BCn-sum}.   This proof is based on $q$-difference equations for the Jackson integrals 
arising from the the two-term relations of Section \ref{Sect.4}.  
\par
We remark that our fundamental invariants $E_r(z)=E_r(a_1,a_2;z)$ %($r=0,1,\ldots,n$) 
are essentially the interpolation theta functions of Coskun--Gustafson \cite{CG2007} and Rains \cite{Ra2006} attached to single columns. In fact, we have 
\begin{align*}
E_r(a_1,a_2;z)&=c_r\, W_{(1^r)}(z/a_2;q,p,t,a_2^2,a_2/qa_1)
{\prod_{i=1}^n \theta(a_1z_i;p)\theta(a_1/z_i;p)} \\
&=c^*_r\, P_{(1^r)}^{*(1,n)}(z;a_2,a_1,q,t;p) \qquad (r=0,1,\ldots,n),
\end{align*}
with the notations $W_\lambda(x;q,p,t,a,b)$ and $P^{*(m,n)}_\lambda(x;a,b,q,t;p)$ of \cite{CG2007} and \cite{Ra2006}, respectively,  
where the constants $c_r$ and $c^*_r$ are determined by comparing the values at 
$
z=(qa_2t^{n-1},\ldots,qa_2t^{n-r}, a_2t^{n-r-1},\ldots,a_2)$.% 
\footnote{
The values of $W_{(1^r)}(x;q,p,t,a,b)$ and $P^{*(1,n)}_{(1^r)}(x;a,b,q,t;p)$ at this point are explicitly given in \cite{CG2007,Ra2006}, while 
the value of 
$
E_r(a_1,a_2;z)
$ is computed by 
$$
E_r(a_1,a_2;z_1,\ldots,z_r,a_2t^{n-r-1},\ldots,a_2)= 
\prod_{i=1}^r
\frac{\theta(a_2t^{n-r}z_i;p)\theta(a_2t^{n-r}/z_i;p)}
{\theta(a_1a_2t^{n-i};p)\theta(a_2t^{n-2r+i}/a_1;p)}.
$$
}
\par
We expect that our method is applicable to other elliptic Jackson summations as well. 
Interesting cases to consider would be the elliptic Gustafson--Rakha summation on $A_m$ in \cite{Sp2003} 
and its variations in \cite{SW2011}, which have remained unproven until now.
%
%
%One of the interesting cases to consider would be the elliptic analogue of 
%the Gustafson--Rakha summation on $A_m$ \cite{SW2011}, which has remained unproven until now. 

\section{Jackson integrals}
\label{Sect.2}

We first rewrite the summation formula \eqref{BCn-sum} in the terminology 
of Jackson integrals. 
For a function $f(z)$ in $n$ variables $z=(z_1,\ldots,z_n)$, we denote by
\begin{equation*}
\int_{0}^{\xi\infty} f(z)\,\varpi_q(z)
=\sum_{\nu\in \mathbb{Z}^n}{}f(q^\nu \xi)
=
\sum_{(\nu_1,\ldots,\nu_n)\in\mathbb{Z}^n}\ f(q^{\nu_1}\xi_1,\ldots,
q^{\nu_n}\xi_n)
\end{equation*}
the Jackson integral of $f(z)$ 
with a base point $\xi=(\xi_1,\ldots,\xi_n)\in(\mathbb{C}^\ast)^n$, 
where 
\begin{equation*}
\varpi_q(z)=\frac{1}{(1-q)^n}\frac{d_qz_1}{z_1}\wedge\cdots\wedge
\frac{d_qz_n}{z_n}. 
\end{equation*}
\par
We define a multivalued meromorphic function 
$\Phi(z)=\Phi(a_1,\ldots,a_6;z)$ on the algebraic torus 
$(\mathbb{C}^\ast)^n$ 
by
\begin{eqnarray}\label{def-Phi}
&&
\Phi(z)=\prod_{i=1}^{n}
\prod_{m=1}^{6}
z_i^{\frac{1}{2}-\alpha_m} \frac{\Gamma(a_mz_i;p,q)}
{\Gamma(qa_m^{-1}z_i;p,q)}
\prod_{1\le j<k\le n}
\frac
{z_j^{1-2\tau}\Gamma(tz_j/z_k;p,q)\Gamma(t z_jz_k;p,q)}
{\Gamma(qt^{-1}z_j/z_k;p,q)\Gamma(qt^{-1} z_jz_k;p,q)}
\nonumber\\
&&\hspace{40pt}
\times
\prod_{i=1}^{n}
z_i^{-1}\theta(z_i^2;p)
\prod_{1\le j<k\le n}
z_j^{-1}\theta(z_j/z_k;p)\theta(z_jz_k;p)
\hspace{10pt}
\end{eqnarray}
with complex parameters $a_i=q^{\alpha_i}$ ($i=1,\ldots,6$)
and $t=q^\tau$, where 
\begin{equation*}
\Gamma(x;p,q)=\frac{(pq/x;p,q)_\infty}{(x;p,q)_{\infty}},
\quad (x;p,q)_\infty=\prod_{i,j=0}^{\infty}(1-p^iq^jx)\quad(|p|<1, |q|<1)
\end{equation*}
denotes the Ruijsenaars elliptic gamma function \cite{Ru1997}.  
Note that
\begin{equation*}
\frac{\Gamma(qx;p,q)}{\Gamma(x;p,q)}=\theta(x;p),\quad
\frac{\Gamma(q^k x;p,q)}{\Gamma(x;p,q)}=\theta(x;p;q)_k
\quad(k=0,1,2,\ldots). 
\end{equation*}
We remark that 
for each $\nu\in\mathbb{Z}^n$ with $\nu_1\ge \nu_2\ge\cdots\ge \nu_n$,
\begin{eqnarray*}
&&\frac{\Phi(q^{\nu}z)}{\Phi(z)}
=
\prod_{i=1}^{n}
\frac{q^{2\nu_i}\theta(q^{2\nu_i}z_i^2;p)}{\theta(z_i^2;p)}
\prod_{1\le j<k\le n}\,
\frac
{\theta(q^{\nu_j-\nu_k}z_j/z_k;p)\theta(q^{\nu_j+\nu_k}z_jz_k;p)}
{\theta(z_j/z_k;p)\theta(z_jz_k;p)}
\nonumber\\
&&
\hspace{15pt}
\times
\prod_{i=1}^{n}
\prod_{m=1}^{6}
\hspace{-1pt}
\frac{a_m^{-\nu_i}\,
\theta(a_mz_i;p;q)_{\nu_i}}
{\theta(qa_m^{-1}z_i;p;q)_{\nu_i}
}
%\hspace{-2pt}
\prod_{1\le j<k\le n}
%\hspace{-1pt}
\frac
{t^{-2\nu_j}\,\theta(tz_j/z_k;p;q)_{\nu_j-\nu_k}\theta(t z_jz_k;p;q)_{\nu_j+\nu_k}%\hspace{8pt}
}
{\theta(qt^{-1}z_j/z_k;p;q)_{\nu_j-\nu_k}
\theta(qt^{-1}z_jz_k;p;q)_{\nu_j+\nu_k}},  
\end{eqnarray*}
and hence 
the summand of the left-hand side of \eqref{BCn-sum} 
coincides with the value of $\Phi(q^\nu z)/\Phi(z)$ at $z=\xi$, 
under the condition $a_1\cdots a_6\,t^{2n-2}=q$. 
\par
In order to rewrite \eqref{BCn-sum} as a Jackson integral, we consider the 
multiplicative lattice
\begin{equation*}
\Lambda_{\xi}=
\left\{\ \xi q^{\nu}=(\xi_1q^{\nu_1},\xi_2 q^{\nu_2},\ldots,
\xi_n q^{\nu_n})\ |\ 
\nu=(\nu_1,\nu_2,\ldots,\nu_n)\in\mathbb{Z}^n\ \right\}
\end{equation*}
with the base point $\xi=(a_1t^{n-1},a_1t^{n-2},\ldots,a_1t,a_1)$, 
under the condition $a_1a_6t^{n-1}=q^{-N}$ 
which makes the Jackson integral a finite sum. 
Notice that $\Phi(z)$ 
has poles along $a_6z_1=q^{k}$ 
($k\in\mathbb{Z}_{\le 0}$) which 
arise from the factor
\begin{equation*}
\Gamma(a_6 z_1;p,q)=
\frac{(pq/a_6z_1;p,q)_\infty}{(a_6z_1;p,q)_\infty}. 
\end{equation*}
Hence, 
the points $z=(a_1t^{n-1}q^{\nu_1},\ldots,a_1q^{\nu_n})\in\Lambda_\xi$ 
with $\nu_1\le N$ are contained in the set of the poles of $\Phi(z)$, 
since 
$a_6z_1=q^{-N+\nu_1}$ 
for $z_1=a_1t^{n-1}q^{\nu_1}$. 
So as to avoid these poles, we modify $\Phi(z)$ as 
\begin{equation*}
\widetilde{\Phi}(z)
=z_1^{-1+2\alpha_6}\dfrac{\theta(a_6z_1;q)}{\theta(qa_6^{-1}z_1;q)}
\Phi(z)
\end{equation*}
by multiplying with a $q$-periodic function. 
Then we see that 
\begin{equation}\label{wPhi-Phi}
\frac{\widetilde{\Phi}(q^{\nu}z)}{\widetilde{\Phi}(z)}=
\frac{\Phi(q^{\nu}z)}{\Phi(z)}\qquad(\nu\in\mathbb{Z}^n)
\end{equation}
and that the lattice $\Lambda_\xi$ 
is contained in the domain of holomorphy 
of $\widetilde{\Phi}(z)$.
Accordingly, 
it makes sense to consider the Jackson integral 
\begin{equation}\label{Jackson-sum}
\dfrac{1}{\widetilde{\Phi}(\xi)}\int_0^{\xi\infty}\widetilde{\Phi}(z)\varpi_q(z),
\qquad \xi=(a_1t^{n-1},a_1t^{n-2},\ldots,a_1)
\end{equation}
under the condition  $a_1a_6t^{n-1}=q^{-N}$. 
Furthermore, 
one can verify that $\widetilde{\Phi}(z)$ has zeros along 
\begin{eqnarray*}
&&z_n=a_1 q^{k}\quad (k<0),\quad
z_j/z_{j+1}=tq^{k}\quad(1\le j\le n-1;\ k<0), 
\nonumber\\
&&
a_6z_1
=q^{k}\quad(k>0), 
\end{eqnarray*}
which implies 
$\widetilde{\Phi}(q^{\nu}\xi)=0$ unless $N\ge \nu_1\ge\nu_2\ge\cdots\ge \nu_n\ge 0$.  
Hence we see that the Jackson integral \eqref{Jackson-sum} is a finite sum, and 
coincides with 
the left-hand side of \eqref{BCn-sum} under the condition
$a_1\cdots a_6\,t^{2n-2}=q$. 
Note also that, if $N=0$, i.e. $a_1a_6t^{n-1}=1$, this sum reduces to a single term, giving 1. 

\section{Fundamental \boldmath$BC_n$ invariants}
\label{Sect.3}

In this section, we introduce a family of fundamental $BC_n$ invariant functions $E_r(a_1,a_2;z)$ $(r=0,\ldots, n)$, 
which are characterized by an interpolation property.
\par\medskip
We denote by $W_n=(\mathbb{Z}/2\mathbb{Z})^n\rtimes \mathfrak{S}_n$ the Weyl group of type $C_n$. 
This group acts on the space ${\cal O}(({\Bbb C}^*)^n)$ of holomorphic functions on $({\Bbb C}^*)^n$ 
through permutations and inversions of the variables $z_1,\ldots,z_n$. 
We consider the $\mathbb{C}$-subspace 
$H_n\subset  {\cal O}(({\Bbb C}^*)^n)$ consisting of all $W_n$-invariant holomorphic functions $f(z)$ such that $T_{p,z_i}f(z)=f(z)/(pz_i^2)\ (i=1,\ldots,n)$, where $T_{p,z_i}$ stands for the $p$-shift operator 
in $z_i$ defined by $T_{p,z_i}f(z_1,\ldots, z_n)=f(\ldots,pz_i,\ldots)$:
\begin{equation}\label{eq:def-Hn}
H_n=\left\{\ f(z)\in\mathcal{O}((\mathbb{C}^\ast)^n)^{W_n}\mid  
T_{p,z_i}f(z)=f(z)/(pz_i^2)\ \ (i=1,\ldots,n)
\right\}. 
\end{equation}
The functions in $H_n$ are called {\it the $BC_n$-symmetric theta functions of degree 1} \cite{Ra2010} and 
the dimension of $H_n$ as a $\mathbb{C}$-vector space is known to be $n+1$, see e.g.  \cite{Ito2008,Ra2010}.
\par
For generic complex parameters $(a_1, a_2)$ and $t$, 
we define a set of reference points 
$\zeta^{(s)}=\zeta^{(s)}(a_1,a_2)\in(\mathbb{C}^\ast)^n$ 
($s=0,1,\ldots,n$) by
\begin{equation*}
\zeta^{(s)}=(\,
\underbrace{a_1,a_1t,\ldots,a_1t^{s-1}}_{s},
\underbrace{ a_2,a_2t,\ldots,a_2t^{n-s-1}}_{n-s}\,)\in (\mathbb{C}^\ast)^n.
\end{equation*}
\begin{thm} \label{thm:E2-basis}
The $\mathbb{C}$-linear space $H_n$ 
has a unique basis $E_r(z)=E_r(a_1,a_2;z)$
$(r=0,1,\ldots,n)$ 
such that 
\begin{equation}
\label{eq:E2-vanish} 
E_r(\zeta^{(s)})
=\delta_{rs} \quad\mbox{for}\quad 0\le r,s\le n,
\end{equation}
where $\delta_{rs}$ is the Kronecker delta.  
\end{thm}
We call the holomorphic functions 
$E_r(z)$ ($r=0,1,\ldots, n$) of 
Theorem \ref{thm:E2-basis}
the {\it $($elliptic$)$ fundamental $W_n$-invariants}. 
We prove this theorem by explicitly constructing 
$E_r(z)$ ($r=0,1,\ldots,n$) with the interpolation property
\eqref{eq:E2-vanish}. 
In what follows, we use the abbreviation $\theta(ax^{\pm};p)=
\theta(ax;p)\theta(ax^{-1};p)$ for the product of two factors with $x$ and $x^{-1}$. 
\begin{df}
{\rm
For each $r=0,1,\ldots,n$, we define a holomorphic function 
$E_r(z)=E_r(z_1,\ldots,z_n)$ 
in $n$ variables $z=(z_1,\ldots,z_n)\in(\mathbb{C}^\ast)^n$ 
with parameters $(a_1,a_2,t)$ by 
\begin{equation*}
E_{r}(z)=
\sum_{\substack{I\subseteq\{1,\ldots,n\}\\ |I|=r}}
\prod_{k=1}^{r}\frac{\theta(a_2 t^{i_k-k}z_{i_k}^{\pm 1};p)}{\theta(a_2 t^{i_k-k}(a_1t^{k-1})^{\pm 1};p)}
\prod_{l=1}^{n-r}\frac{\theta(a_1t^{j_l-l}z_{j_l}^{\pm 1};p)}{\theta(a_1t^{j_l-l}(a_2t^{l-1})^{\pm 1};p)},
\end{equation*}
where  the summation is over all $r$-subsets $I$ of $\{1,\ldots,n\}$, and  
$I=\{i_1<\cdots<i_r\}$, $J=\{1,\ldots,n\}\backslash I=\{ j_1<\cdots<j_{n-r}\}$.
}
\end{df}
In particular, 
\begin{equation*}
E_0(z)=
\prod_{i=1}^{n}\frac{\theta(a_1z_{i}^{\pm 1};p)}
{\theta(a_1(a_2t^{i-1})^{\pm 1};p)},
\quad
E_n(z)=
\prod_{i=1}^{n}\frac{\theta(a_2z_{i}^{\pm 1};p)}
{\theta(a_2(a_1t^{i-1})^{\pm 1};p)}. 
\end{equation*}

When we need to specify the number of variables, 
we use the notation $E^{(n)}_r(z)=E_r(z)$.  
We first remark that these functions can also be 
defined inductively on the number of variables. 
\begin{prop}
\label{prop:E2-rec} 
The functions $E^{(n)}_r(z)$ $(r=0,1,\ldots,n)$
satisfy the recurrence relation 
\begin{equation*}
E^{(n)}_{r}(z)%&
=E^{(n-1)}_{r-1}(z_{\widehat{n}})
\frac{\theta(a_2t^{n-r}z_n^{\pm 1};p)}
{\theta(a_2t^{n-r}(a_1t^{r-1})^{\pm 1};p)}
+E^{(n-1)}_{r}(z_{\widehat{n}})
\frac{\theta(a_1t^{r}z_n^{\pm 1};p)}
{\theta(a_1t^{r}(a_2t^{n-r-1})^{\pm 1};p)},
\end{equation*}
where $z_{\,\widehat{i}}=(z_1,\ldots,z_{i-1},z_{i+1},\ldots,z_n)\quad (i=1,\ldots,n)$. 
\hfill $\square$
\end{prop}
\par
It is immediately verified that these functions 
have the quasi-periodicity 
\begin{equation*}
T_{p,z_i}E_r(z)=E_r(z)/(pz_i^2)\quad (i=1,\ldots,n).  
\end{equation*}
We next prove that $E_r(z)$ ($r=0,1,\ldots,n$) 
carry the interpolation property 
with respect to the reference points $\zeta^{(s)}$ ($s=0,1,\ldots,n$), 
\begin{lem}
\label{lem:E2-vanish} \ \ 
$E_r(\zeta^{(s)})=\delta_{rs}$ 
for $0\le r,s\le n$.  
\end{lem}
{\bf Proof.}
We take   
the summand of $E_r(z)$ corresponding to $I$, 
and look at its value at $\zeta^{(s)}$.  This term 
has zeros along the divisors
\begin{equation*}
z_{i_k}=a_2t^{i_k-k}\quad(k=1,\ldots,r),\quad
z_{j_l}=a_2t^{j_l-l}\quad(l=1,\ldots,n-r).  
\end{equation*}
On the other hand, the components of $\zeta^{(s)}$ are given by 
\begin{equation*}
\zeta^{(s)}_i=a_1t^{i-1}\quad(1\le i\le s),\quad
\zeta^{(s)}_i=a_2t^{i-s-1}\quad(s+1\le i\le n). 
\end{equation*}
When $r>s$, we have 
$\zeta^{(s)}_{i_{s+1}}=a_2t^{i_{s+1}-s-1}$ since $i_{s+1}\ge s+1$, 
and hence the factor for $k=s+1$ vanishes. 
When $r<s$,  we have 
$\zeta^{(s)}_{j_1}=a_1t^{j_1-1}$ since $j_1\le r+1\le s$, 
and hence the factor for $l=1$ vanishes.  
When $r=s$ and $j_1\le r$, we have 
$\zeta^{(r)}_{j_1}=a_1t^{j_1-1}$, and the factor for $l=1$ 
vanishes.  
Hence the summand 
is nontrivial only when $r=s$ and $j_1=r+1$, namely 
$\{i_1,\ldots,i_r\}=\{1,\ldots,r\}$ and $\{j_1,\ldots,j_{n-r}\}=\{r+1,\ldots,n\}$; 
in this case the corresponding term takes the value 1 since 
\begin{equation*}
\hspace{30pt}
\zeta^{(r)}_k=a_2t^{k-1}\quad(k=1,\ldots,r),
\quad
\zeta^{(r)}_{r+l}=a_1t^{l-1}\quad(l=1,\ldots,n-r).
\hspace{30pt}\square
\end{equation*}
\par\medskip
Notice that 
the linear independence of $E_r(z)$ ($r=0,1,\ldots,n$) 
follows from the interpolation property.  
In order to establish Theorem \ref{thm:E2-basis}, 
it remains to show the $W_n$-invariance of $E_r(z)$. 
Since $E_r(z)$ are invariant under the inversion of 
variables $z_i$, we need to prove they are symmetric 
in the variables $(z_1,\ldots,z_n)$.
\begin{lem}
\label{lem:E2-Winv}
For each $r=0,1,\ldots,n$, $E_r(z)$ is symmetric in the 
variables $z_1,\ldots,z_n$. 
\end{lem}
{\bf Proof.}
We prove this lemma by induction on the number of 
variables $n$.   
By the recurrence relation of Proposition \ref{prop:E2-rec}, 
$E^{(n)}_r(z)=E^{(n)}_r(z_{\widehat{n}},z_n)$ are symmetric 
with respect to $z_{\widehat{n}}=(z_1,\ldots,z_{n-1})$.  
Hence it suffices to show that they are invariant under the 
transposition of $(z_{n-1}, z_n)$, namely, 
$E^{(n)}_r(z_{\widehat{n-1,n}},z_{n-1},z_n)=E^{(n)}_r(z_{\widehat{n-1,n}},z_n,z_{n-1})$,
where  
$z_{\widehat{n-1,n}}=(z_1,\ldots,z_{n-2})$.
Using Proposition \ref{prop:E2-rec} twice, 
we have
\begin{equation*}
\begin{split}
E^{(n)}_r(z)
&=
E^{(n-2)}_{r-2}(z_{\widehat{n-1,n}})
\dfrac{\theta(a_2t^{n-r}(z_{n-1})^{\pm1 };p)}
{\theta(a_2t^{n-r}(a_1t^{r-2})^{\pm1 };p)}
\dfrac{\theta(a_2t^{n-r}(z_{n})^{\pm1 };p)}
{\theta(a_2t^{n-r}(a_1t^{r-1})^{\pm1 };p)}
\nonumber\\
&
\quad\mbox{}+
E^{(n-2)}_{r-1}(z_{\widehat{n-1,n}})
\Bigg(
\frac{\theta(a_1t^{r-1}(z_{n-1})^{\pm1 };p)}
{\theta(a_1t^{r-1}(a_2t^{n-r-1})^{\pm1 };p)}
\frac{\theta(a_2t^{n-r}(z_{n})^{\pm1 };p)}
{\theta(a_2t^{n-r}(a_1t^{r-1})^{\pm1 };p)}
\nonumber\\
&%
\hspace{90pt}
+
\frac{\theta(a_2t^{n-r-1}(z_{n-1})^{\pm1 };p)}
{\theta(a_2t^{n-r-1}(a_1t^{r-1})^{\pm1 };p)}
\frac{\theta(a_1t^{r}(z_{n})^{\pm1 };p)}
{\theta(a_1t^{r}(a_2t^{n-r-1})^{\pm1 };p)}
\Bigg) 
\nonumber\\
&\quad\mbox{}+
E^{(n-2)}_{r}(z_{\widehat{n-1,n}})
\dfrac{\theta(a_1t^{r}(z_{n-1})^{\pm1 };p)}
{\theta(a_1t^{r}(a_2t^{n-r-2})^{\pm1 };p)}
\frac{\theta(a_1t^{r}(z_{n})^{\pm1 };p)}
{\theta(a_1t^{r}(a_2t^{n-r-1})^{\pm1 };p)}
.
\end{split}
\end{equation*}
Hence we obtain
\begin{equation*}
\begin{split}
&E^{(n)}_r(z_{\widehat{n-1,n}},z_{n-1},z_{n})-
E^{(n)}_r(z_{\widehat{n-1,n}},z_{n},z_{n-1})
\\
&\hspace{10pt}
=
\frac{E^{(n-2)}_{r-1}(z_{\widehat{n-1,n}})}{\theta(a_1t^{r-1}(a_2t^{n-r-1})^{\pm1 };p)
\theta(a_2t^{n-r}(a_1t^{r-1})^{\pm1 };p)}
\\
&\hspace{60pt}
\times\bigg(
\theta(a_1t^{r-1}(z_{n-1})^{\pm1 };p)
\theta(a_2t^{n-r}(z_{n})^{\pm1 };p)\\
&\hspace{80pt}
-
\theta(a_1t^{r-1}(z_{n})^{\pm1 };p)
\theta(a_2t^{n-r}(z_{n-1})^{\pm1 };p)\bigg)
\\
&\hspace{20pt}
+
\frac{E^{(n-2)}_{r-1}(z_{\widehat{n-1,n}})}{\theta(a_2t^{n-r-1}(a_1t^{r-1})^{\pm1 };p)
\theta(a_1t^{r}(a_2t^{n-r-1})^{\pm1 };p)}\\
&\hspace{60pt}
\times\bigg(
\theta(a_2t^{n-r-1}(z_{n-1})^{\pm1 };p)
\theta(a_1t^{r}(z_{n})^{\pm1 };p)\\
&\hspace{80pt}
-\theta(a_2t^{n-r-1}(z_{n})^{\pm1 };p)
\theta(a_1t^{r}(z_{n-1})^{\pm1 };p)
\bigg). 
\end{split}
\end{equation*}
By the three-term relation
\begin{equation*}
\theta(xu^{\pm 1};p)\theta(yv^{\pm 1};p)
-\theta(xv^{\pm 1};p)\theta(yu^{\pm 1};p)
=\dfrac{y}{u}\theta(xy^{\pm 1};p)\theta(uv^{\pm 1};p)
\end{equation*}
of the theta function, we obtain
\begin{equation*}
\begin{split}
&E^{(n)}_r(z_{\widehat{n-1,n}},z_{n-1},z_{n})-
E^{(n)}_r(z_{\widehat{n-1,n}},z_{n},z_{n-1})
\nonumber\\
&\hspace{10pt}
=-E^{(n-2)}_{r-1}(z_{\widehat{n-1,n}})
\Bigg(
\dfrac{a_1t^{r-1}
\theta(z_{n-1}(z_{n})^{\pm1 };p)
}{z_{n-1}
\theta(a_1t^{r-1}(a_2t^{n-r-1})^{\pm1 };p)
}\\
&\hspace{105pt}+
\dfrac{a_2t^{n-r-1}
\theta(z_{n-1}(z_{n})^{\pm1 };p)
}{z_{n-1}
\theta(a_2t^{n-r-1}(a_1t^{r-1})^{\pm1 };p)
}
\Bigg)\\
&\hspace{10pt}
=0
\end{split}
\end{equation*}
as desired. 
\hfill$\square$
\par\medskip
We have thus proved that $E_r(z)=E_r(a_1,a_2;z)$ $(r=0,1\ldots,n)$ provide a basis of $H_n$ satisfying the condition \eqref{eq:E2-vanish}.
The uniqueness of such a basis also follows from the interpolation property. 
This completes the proof of Theorem \ref{thm:E2-basis}. 
\section{Two-term relations}
\label{Sect.4}
As we have seen in Section 2, the left-hand side of  
\eqref{BCn-sum} coincides with the Jackson integral 
\eqref{Jackson-sum}.  
Hence our goal is to prove that if the balancing condition 
$a_1\cdots a_6 t^{2n-2}=q$ is satisfied then
\begin{equation}
\label{BCn-Jackson-sum}
\int_0^{\xi\infty}
\dfrac{\widetilde{\Phi}(z)}{\widetilde{\Phi}(\xi)}\varpi_q(z)
=
\prod_{i=1}^{n}
\frac{
\theta(qa_1^2t^{n+i-2};p;q)_N
\prod_{2\le j<k\le 4}\theta(qa_j^{-1}a_k^{-1}t^{1-i};p;q)_N
}{
\theta(qa_1^{-1}a_2^{-1}a_3^{-1}a_4^{-1}t^{2-n-i};p;q)_N
\prod_{m=2}^{4}\theta(qa_m^{-1}a_1t^{i-1};p;q)_N
},
\end{equation}
where $\xi=(a_1t^{n-1},a_1t^{n-2},\ldots,a_1)$ and  
$a_1a_6t^{n-1}=q^{-N}$. 
For this purpose, we investigate the $q$-difference equation 
to be satisfied by the left-hand side.  
In what follows, we always assume  
$\xi=(a_1t^{n-1},a_1t^{n-2},\ldots,a_1)$, $a_1a_6t^{n-1}=q^{-N}$ 
($N=0,1,2,\ldots$)
so that the Jackson integral makes sense as a finite sum.  
However, we do not impose 
any balancing condition on the parameters 
$(a_1,\ldots,a_6)$
until it becomes necessary.  
\par\medskip
Let $\varphi(z)$ be a meromorphic function on 
$(\mathbb{C}^\ast)^n$, and suppose that 
$\varphi(z)$ is holomorphic 
in a neighborhood of the multiplicative lattice 
$\Lambda_{\xi}=\{\, q^{\nu}\xi\mid \nu\in\mathbb{Z}^n\,\}$. 
Then by the $q$-shift invariance of Jackson integrals, 
for $i=1,\ldots,n$ we have 
\begin{align*}
\int_0^{\xi\infty}\varphi(z)
\dfrac{\widetilde{\Phi}(z)}{\widetilde{\Phi}(\xi)}\varpi_q(z)
&=
\int_0^{\xi\infty}T_{q,z_i}
\left(
\varphi(z)
\dfrac{\widetilde{\Phi}(z)}{\widetilde{\Phi}(\xi)}
\right)
\varpi_q(z)
\nonumber\\
&=
\int_0^{\xi\infty}T_{q,z_i}\varphi(z)
\dfrac{T_{q,z_i}\widetilde{\Phi}(z)}
{\widetilde{\Phi}(\xi)}
\varpi_q(z), 
\end{align*}
and hence
\begin{equation*}
\int_0^{\xi\infty}
\left(
\varphi(z)-
\dfrac{T_{q,z_i}\widetilde{\Phi}(z)}{\widetilde{\Phi}(z)}
T_{q,z_i}\varphi(z)
\right)
\dfrac{\widetilde{\Phi}(z)}{\widetilde{\Phi}(\xi)}\varpi_q(z)
=0.  
\end{equation*}
In view of this property, we introduce 
the operator $\nabla_{q,z_i}$ 
by setting 
\begin{equation}\label{eq:nabla}
\nabla_{q,z_i}\varphi(z)=
\varphi(z)-
\dfrac{T_{q,z_i}\widetilde{\Phi}(z)}{\widetilde{\Phi}(z)}
T_{q,z_i}\varphi(z)
\end{equation}
for $i=1,\ldots,n$. 
In the notation 
\begin{equation*}
\langle \varphi(z)\rangle = 
\int_0^{\xi\infty}\varphi(z)
\dfrac{\widetilde{\Phi}(z)}{\widetilde{\Phi}(\xi)}\varpi_q(z)
\end{equation*}
of expectation values, the $q$-shift invariance of 
Jackson integrals is expressed as 
\begin{equation}
\label{eq:q-Stokes}
\langle \nabla_{q,z_i}\varphi(z)\rangle=0\quad(i=1,\ldots,n)
\end{equation}
provided that 
$\varphi(z)$ is holomorphic 
in a neighborhood of $\Lambda_\xi$. 
In our case, 
from \eqref{def-Phi} and \eqref{wPhi-Phi}
the ratio 
$T_{q,z_i}\widetilde{\Phi}(z)/\widetilde{\Phi}(z)=T_{q,z_i}\Phi(z)/\Phi(z)$ 
is explicitly written as 
\begin{equation}
\label{eq:ratio}
\frac{T_{q,z_i}\widetilde{\Phi}(z)}{\widetilde{\Phi}(z)}
=
\frac{\theta(q^2z_i^2;p)}{q\,\theta(z_i^2;p)}
\prod_{m=1}^{6}\frac{q^{\frac{1}{2}}a_m^{-1}\theta(a_mz_i;p)}
{\theta(qa_m^{-1}z_i;p)}
\prod_{\substack{1\le k\le n\\
k\ne i}}
\frac{\theta(tz_iz_k^{\pm1};p)
\theta(q^{-1}z_i^{-1}z_k^{\pm1};p)}
{\theta(z_iz_k^{\pm1};p)\theta(q^{-1}tz_i^{-1}z_k^{\pm 1};p)
}. 
\end{equation}
\par
In this section we prove two-term relations for the 
Jackson integrals $\langle E_r(z)\rangle$ 
of the fundamental $W_n$-invariants 
$E_r(z)=E_r(a_1,a_2;z)$ ($r=0,1,\ldots,n$).
\begin{thm}
\label{thm:2-term} 
Suppose that $a_1a_2\cdots a_{6}t^{2n-2}=1$. 
Then
\begin{equation}
\label{eq:2-term}
\la E_{r}(a_1,a_2;z)\ra=-c_{r}\la E_{r-1}(a_1,a_2;z)\ra
\quad (r=1,\ldots,n), 
\end{equation}
where the coefficients $c_r$ are given by 
\begin{equation*}
c_{r}=
\frac{
a_1^2t^{2r-2}\,
\theta(t^{n-r+1};p)
\theta(a_2a_1^{-1}t^{n-r+1};p)
\theta(a_1a_2^{-1}t^{-n+2r};p)
}
{
a_2^2t^{2n-2r}\,
\theta(t^r;p)
\theta(a_2a_1^{-1}t^{n-2r+2};p)
\theta(a_1a_2^{-1}t^r;p)
}
\prod_{m=3}^6
\dfrac{\theta(a_ma_2t^{n-r};p)}{\theta(a_ma_1t^{r-1};p)}.
\end{equation*}
\end{thm}
{\bf Proof.}\ 
We remark that \eqref{eq:ratio} can be rewritten as
\begin{equation*}
\begin{split}
\frac{T_{q,z_i}\widetilde{\Phi}(z)}{\widetilde{\Phi}(z)}
=&
-\frac{(q^{-1}z_i^{-1})^2\theta(q^{-2}z_i^{-2};p)}
{z_i^2\,\theta(z_i^2;p)}
\prod_{m=1}^{6}\frac{
\theta(a_mz_i;p)}{
\theta(q^{-1}a_mz_i^{-1};p)}\\
&\quad\times
\prod_{\substack{1\le k\le n\\
k\ne i}}
\frac{\theta(tz_iz_k^{\pm1};p)
\theta(q^{-1}z_i^{-1}z_k^{\pm1};p)}
{\theta(z_iz_k^{\pm1};p)\theta(q^{-1}tz_i^{-1}z_k^{\pm 1};p)
}. 
\end{split}
\end{equation*}
Hence the above ratio
is expressed as 
\begin{equation}
\label{eq:ratio2}
\frac{T_{q,z_i}\widetilde{\Phi}(z)}{\widetilde{\Phi}(z)}
=-\frac{F_i(z)}{T_{q,z_i}(F_i(z^{-1}))},
\end{equation}
where $z^{-1}=(z_1^{-1},\ldots,z_n^{-1})$ and
\begin{equation}\label{eq:Fi}
F_{i}(z)=
\frac{\prod_{m=1}^{6}
\theta(a_mz_i;p)}{z_i^2\, \theta(z_i^2;p)}
\prod_{\substack{1\le j\le n\\ j\ne i}}
\frac{\theta(tz_iz_j^{\pm 1};p)}{\theta(z_iz_j^{\pm 1};p)}. 
\end{equation}
\par
Fixing an index $r\in\{1,\ldots,n\}$, 
for each $i=1,\ldots,n$ we define a meromorphic function 
$\varphi_i(z)$ by 
\begin{equation*}
\varphi_i(z)=
F_{i}(z^{-1})E_{r-1}^{(n-1)}(z_{\,\widehat{i}\,}), 
\end{equation*}
where $z_{\,\widehat{i}}=
(z_1,\ldots,z_{i-1},z_{i+1},\ldots,z_n)$. 
Note that $\varphi_i(z)$ is holomorphic in a neighborhood 
of $\Lambda_{\xi}$, $\xi=(a_1t^{n-1},a_1t^{n-2},\ldots,a_1)$.  
Then, from the definition \eqref{eq:nabla} of 
$\nabla_{q,z_i}$ and \eqref{eq:ratio2}, we have 
\begin{equation*}
\nabla_{q,z_i}\,\varphi_i(z)=
\Big(F_{i}(z^{-1})+F_{i}(z)\Big)
E_{r-1}^{(n-1)}(z_{\,\widehat{i}\,}),
\end{equation*}
and hence 
\begin{equation*}
\sum_{i=1}^n \nabla_{q,z_i}\,\varphi_i(z)=
\sum_{i=1}^n\Big(F_{i}(z^{-1})+F_{i}(z)\Big)
E_{r-1}^{(n-1)}(z_{\,\widehat{i}\,}). 
\end{equation*}
Here we set
\begin{equation}\label{eq:def-hr}
h_r(z)=
\sum_{i=1}^n\Big(F_{i}(z^{-1})+F_{i}(z)\Big)
E_{r-1}^{(n-1)}(z_{\,\widehat{i}\,}) 
\end{equation}
for $r=1,\ldots,n$.  
Then we have 
\begin{equation}\label{eq:<hr>}
\la h_r(z)\ra =\sum_{i=1}^{n}\la \nabla_{q,z_i}\varphi_i(z)\ra
=0
\end{equation}
by \eqref{eq:q-Stokes}.  
Theorem \ref{thm:2-term} is obtained by combining \eqref{eq:<hr>} with 
Lemma \ref{lem:key} below. 
\hfill$\square$ 
\begin{lem}
\label{lem:key}
Suppose that $a_1a_2\cdots a_{6}t^{2n-2}=1$. 
\newline
$(1)$\ \ For each $r=1,\ldots,n$, 
$h_r(z)$ belongs to the $\mathbb{C}$-linear space $H_n$ of \eqref{eq:def-Hn}.
\newline
$(2)$\ \ In terms of the fundamental $W_n$-invariants $E_s(z)$ $(s=0,1,\ldots,n)$, 
$h_r(z)$ is expressed as 
\begin{equation}
\label{eq:key}
h_r(z)=c_{r,r}E_{r}^{(n)}(z)+c_{r,r-1}E_{r-1}^{(n)}(z),
\end{equation}
where the coefficients $c_{r,r}$ and $c_{r,r-1}$ are given by 
\begin{equation}
\label{eq:crr}
\begin{split}
c_{r,r}&=F_r(\zeta^{(r)})\\
&=
\frac{\theta(t^{r};p)\,\theta(a_1a_2t^{n-1};p)\,\theta(a_1a_2^{-1}t^{r};p)}
{(a_1t^{r-1})^2\,\theta(t;p)\,\theta(a_1a_2^{-1}t^{-n+2r};p)}
\prod_{m=3}^{6}
\theta(a_ma_1t^{r-1};p), 
\\[5pt]
c_{r,r-1}&=F_n(\zeta^{(r-1)})\\
&=
\frac
{\theta(t^{n-r};p)\,\theta(a_1a_2t^{n-1};p)\,\theta(a_2a_1^{-1}t^{n-r+1};p)}
{(a_2t^{n-r})^2\,\theta(t;p)\,\theta(a_2a_1^{-1}t^{n-2r+2};p)}
\prod_{m=3}^{6}
\theta(a_ma_2t^{n-r};p). 
\end{split}
\end{equation}
\end{lem}
{\bf Proof.}\ 
(1)  Under the balancing condition $a_1a_2\cdots a_{6}t^{2n-2}=1$, 
by the explicit form \eqref{eq:Fi} of $F_{i}(z)$ 
it is directly verified that $h_r(z)$ has the required quasi-periodicity. 
By the expression \eqref{eq:def-hr}, $h_r(z)$ is symmetric in the variables 
$z_1,\ldots,z_n$.  The invariance with respect to the inversion $z_k\to z_k^{-1}$ 
follows from the $W_{n-1}$-invariance of $E_{r-1}^{(n-1)}(z_{\,\widehat{i}\,})$ 
($i\ne k$) and $F_k(z^{-1})=F_{k}(z)\big\vert_{z_k\to z_k^{-1}}$.  
It remains to be shown that 
$h_r(z)$ is holomorphic on $(\mathbb{C}^\ast)^n$.
By the quasi-periodicity and the $W_n$-invariance, we have only to show that $h_r(z)$ has 
no poles along the divisors $z_1/z_2=1$ and $z_1^2=1$.   It can be confirmed by 
\begin{equation*}
\theta(z_1/z_2;p)\big(F_{1}(z^{\epsilon})E_{r-1}^{(n-1)}(z_{\,\widehat{1}})
+F_{2}(z^{\epsilon})E_{r-1}^{(n-1)}(z_{\,\widehat{2}})\big)\big\vert_{z_1=z_2}=0\quad(\epsilon=\pm1)
\end{equation*}
and by the fact that $z_1^{-1}\theta(z_1^2;p)\left(F_1(z^{-1})+F_1(z)\right)$ is alternating 
with respect to the inversion $z_1\to z_1^{-1}$, respectively. 
\newline
(2)  
By Theorem \ref{thm:E2-basis}, $h_r(z)$ can be expressed as a linear combination of the fundamental 
$W_n$-invariants:
\begin{equation*}
h_r(z)=\sum_{s=0}^{n}\ c_{r,s}\,E_{s}^{(n)}(z),\quad \mbox{where}\quad
c_{r,s}=h_r(\zeta^{(s)})\quad(s=0,1,\ldots,n). 
\end{equation*}
By the definition \eqref{eq:Fi}, $F_i(z^{-1})$ and $F_i(z)$ $(i=1,\ldots,n)$
satisfy the vanishing property 
\begin{equation*}
F_i((\zeta^{(s)})^{-1})=0 \quad(0\le s\le n),\quad
F_i(\zeta^{(s)})=0\quad(0\le s\le n; i\ne s,n).
\end{equation*}
Hence,
\begin{equation*}
\begin{split}
c_{r,0}=h_r(\zeta^{(0)})&=F_n(\zeta^{(0)})E_{r-1}^{(n-1)}(z_{\widehat{n}})\vert_{z=\zeta^{(0)}}=\delta_{r,1} F_n(\zeta^{(0)}),
\nonumber\\
c_{r,n}=h_r(\zeta^{(n)})&=F_n(\zeta^{(n)})E_{r-1}^{(n-1)}(z_{\widehat{n}})\vert_{z=\zeta^{(n-1)}}=\delta_{r,n} F_n(\zeta^{(n)}). 
\end{split}
\end{equation*}
For $s=1,\ldots,n-1$, we have
\begin{equation*}
\begin{split}
c_{r,s}=h_r(\zeta^{(s)})&=F_s(\zeta^{(s)})E_{r-1}^{(n-1)}(z_{\widehat{n}})\vert_{z=\zeta^{(s-1)}}
+F_n(\zeta^{(s)})E_{r-1}^{(n-1)}(z_{\widehat{n}})\vert_{z=\zeta^{(s)}}
\nonumber\\
&=
\delta_{r,s} F_s(\zeta^{(s)})+\delta_{r-1,s}F_n(\zeta^{(s)}). 
\end{split}
\end{equation*}
Namely, 
\begin{equation*}
c_{r,s}=
\delta_{r,s} F_r(\zeta^{(r)})+\delta_{r-1,s}F_n(\zeta^{(r-1)}). 
\end{equation*}
Note that the formula above is valid in the cases $k=0$ and $n$ as well.  
Finally we obtain
\begin{equation*}
h_r(z)=c_{r,r-1} E_{r-1}(z)+c_{r,r} E_{r}(z);\quad
c_{r,r-1}=F_n(\zeta^{(r-1)}),\ \ c_{r,r}=F_r(\zeta^{(r)})
\end{equation*}
for $r=1,\ldots,n$, which give the expressions of \eqref{eq:key} and \eqref{eq:crr}. 
\hfill$\square$ 
\par\medskip
From Lemma \ref{lem:key}, by \eqref{eq:<hr>} we obtain 
$c_{r,r}\la E_r(z)\ra+c_{r,r-1}\la E_{r-1}(z)\ra=\la h_r(z)\ra=0$, 
namely,
$\la E_r(z)\ra=-c_r\la E_{r-1}(z)\ra$ with 
$c_r=c_{r,r}/c_{r,r-1}$ $(r=1,\ldots,n)$.  
This gives the two-term relation \eqref{eq:2-term}. 

Applying these two-term relations repeatedly, we obtain the formula 
connecting $\la E_n(z)\ra$ and $\la E_0(z)\ra$.  
\begin{cor}\label{cor:EnE0}
Suppose that $a_1a_2\cdots a_6t^{2n-2}=1$.  
Then one has 
\begin{equation*}
\la E_n(a_1,a_2;z)\ra
=
\la E_0(a_1,a_2;z)\ra
\prod_{i=1}^{n}
\left(
\frac{a_1^3}{a_2^3}
\frac{
\theta(a_2a_1^{-1}t^{i-1};p)
}
{
\theta(a_1 a_2^{-1}t^{i-1};p)
}
\prod_{m=3}^6
\dfrac{\theta(a_ma_2t^{i-1};p)}{\theta(a_ma_1t^{i-1};p)}
\right).
\end{equation*}
\end{cor}
\section{Proof of the summation formula}
\label{Sect.5}
In this section, we prove the summation formula 
\eqref{BCn-Jackson-sum} on the basis of Theorem 
\ref{thm:2-term} of the previous section.  
\par\medskip
In order to investigate the dependence of the Jackson integral 
on the parameters $(a_5,a_6)$, we set
\begin{equation*}
J(a_5,a_6)=\la 1\ra=
\int_0^{\xi\infty}
\frac{\widetilde{\Phi}(z)}
{\widetilde{\Phi}(\xi)}\varpi_q(z),
\quad\xi=(a_1,a_1t,\ldots,a_1t^{n-1}), 
\end{equation*}
assuming that $a_1a_6t^{n-1}=q^{-N}$, where  $N$ is a nonnegative integer. 
Since 
\begin{equation*}
T_{q,a_k}\widetilde{\Phi}(z)=\widetilde{\Phi}(z)
\,\prod_{i=1}^{n}(-a_k^{-1})\theta(a_kz_i^{\pm 1};p), 
\end{equation*}
we have 
\begin{equation*}
T_{q,a_k}\left(\frac{\widetilde{\Phi}(z)}
{\widetilde{\Phi}(\xi)}\right)
=\frac{\widetilde{\Phi}(z)}{\widetilde{\Phi}(\xi)}
\,\prod_{i=1}^{n}
\frac{\theta(a_kz_i^{\pm 1};p)}{\theta(a_k(a_1t^{i-1})^{\pm 1};p)}
\end{equation*}
for $k=2,\ldots,6$. 
This implies
\begin{equation*}
\begin{split}
J(a_5,qa_6)
&=\Big\la\,\prod_{i=1}^{n}
\frac{\theta(a_6z_i^{\pm 1};p)}
{\theta(a_6(a_1t^{i-1})^{\pm 1};p)}
\Big\ra
=
\la E_n(a_5,a_6;z)\ra
\prod_{i=1}^{n}
\frac{\theta(a_6(a_5t^{i-1})^{\pm 1};p)}
{\theta(a_6(a_1t^{i-1})^{\pm 1};p)},
\\
J(qa_5,a_6)
&=\Big\la\,\prod_{i=1}^{n}
\frac{\theta(a_5z_i^{\pm 1};p)}
{\theta(a_5(a_1t^{i-1})^{\pm 1};p)}
\Big\ra
=
\la E_0(a_5,a_6;z)\ra
\prod_{i=1}^{n}
\frac{\theta(a_5(a_6t^{i-1})^{\pm 1};p)}
{\theta(a_5(a_1t^{i-1})^{\pm 1};p)}.
\end{split}
\end{equation*}
Then by Corollary \ref{cor:EnE0} for $(a_5,a_6)$, 
we obtain
\begin{align}
\label{eq:J=***J}
J(a_5,qa_6)
&=
J(qa_5,a_6)
\prod_{i=1}^{n}
\left(
\frac{a_5^2}{a_6^2}
\frac{\theta(a_5^{-1}a_1t^{i-1};p)}
{\theta(a_6^{-1}a_1t^{i-1};p)}
\prod_{m=2}^{4}
\dfrac{\theta(a_ma_6t^{i-1};p)}{\theta(a_ma_5t^{i-1};p)}
\right)
\nonumber\\
&
=J(qa_5,a_6)
\prod_{i=1}^{n}
\left(
\frac{\theta(a_5a_1^{-1}t^{-i+1};p)}
{\theta(a_6^{-1}a_1t^{i-1};p)}
\prod_{m=2}^{4}
\dfrac
{\theta(a_m^{-1}a_6^{-1}t^{-i+1};p)}
{\theta(a_ma_5t^{i-1};p)}
\right),
\end{align}
where the second equality follows from $\theta(x;p)=-x\theta(x^{-1};p)$ 
and $a_1\ldots a_6t^{2n-2}=1$.  
\par
We now suppose that  $a_1\cdots a_6 t^{2n-2}=q$ and  $a_1a_6t^{n-1}=q^{-N}$.  
Then, 
\begin{equation*}
J(a_5,a_6)
=
J(q^{-1}a_5,qa_6)
\prod_{i=1}^{n}
\left(
\frac
{\theta(a_6^{-1}a_1t^{i-1};p)}
{\theta(q^{-1}a_5a_1^{-1}t^{-i+1};p)}
\prod_{m=2}^{4}
\dfrac
{\theta(q^{-1}a_ma_5t^{i-1};p)}
{\theta(a_m^{-1}a_6^{-1}t^{-i+1};p)}
\right),  
\end{equation*}
which is obtained from \eqref{eq:J=***J} by replacing $a_5$ with $q^{-1}a_5$. 
Using this formula $N$ times we have
\begin{equation*}
\begin{split}
J(a_5,a_6)
=
J(q^{-N}a_5,q^Na_6)
\prod_{i=1}^{n}
\Bigg(
&\frac
{\theta(q^{1-N}a_6^{-1}a_1t^{i-1};p;q)_N}
{\theta(q^{-N}a_5a_1^{-1}t^{-i+1};p;q)_N}\\
&\hspace{8pt}
\times 
\prod_{m=2}^{4}
\frac
{\theta(q^{-N}a_ma_5t^{i-1};p;q)_N}
{\theta(q^{1-N}a_m^{-1}a_6^{-1}t^{-i+1};p;q)_N}
\Bigg)
.  
\end{split}
\end{equation*}
Since 
$a_1(q^{N}a_6)t^{n-1}=1$, 
the Jackson integral $J(q^{-N}a_5,q^{N}a_6)$ reduces to 1
as we mentioned at the end of Section \ref{Sect.2}.  
Hence we obtain
\begin{equation*}
\begin{split}
J(a_5,a_6)
&=
\prod_{i=1}^{n}
\left(
\frac
{\theta(q^{1-N}a_6^{-1}a_1t^{i-1};p;q)_N}
{\theta(q^{-N}a_5a_1^{-1}t^{-i+1};p;q)_N}
\prod_{m=2}^{4}
\dfrac
{\theta(q^{-N}a_ma_5t^{i-1};p;q)_N}
{\theta(q^{1-N}a_m^{-1}a_6^{-1}t^{-i+1};p;q)_N}
\right)
\\
&=
\prod_{i=1}^{n}
\frac
{\theta(qa_1^2t^{n+i-2};p;q)_N
\prod_{2\le j<k\le 4}
\theta(qa_j^{-1}a_k^{-1}t^{-n+i};p;q)_N
}
{\theta(a_5a_6t^{n-i};p;q)_N
\prod_{m=2}^{4}
\theta(qa_m^{-1}a_1t^{n-i};p;q)_N
}, 
\end{split}
\end{equation*}
which coincides with \eqref{BCn-Jackson-sum}. 
This completes the proof of the summation formula \eqref{BCn-sum}. 
\section*{Acknowledgements} 
This work was partially supported by JSPS Kakenhi Grants (S)24224001, (C)25400118 
and (B)15H03626. 
The authors would like to thank the anonymous referees who kindly 
provided valuable comments and suggestions for the 
improvement of our manuscript.  
%\end{acknowledgements}

% BibTeX users please use one of
%\bibliographystyle{spbasic}      % basic style, author-year citations
%\bibliographystyle{spmpsci}      % mathematics and physical sciences
%\bibliographystyle{spphys}       % APS-like style for physics
%\bibliography{}   % name your BibTeX data base

\begin{thebibliography}{}
%
% and use \bibitem to create references. Consult the Instructions
% for authors for reference list style.
%

%\bibitem{RefJ}
% Format for Journal Reference
%Author, Article title, Journal, Volume, page numbers (year)
% Format for books
%\bibitem{RefB}
%Author, Book title, page numbers. Publisher, place (year)
% etc

%
\bibitem{vDS2000}
J.~F.~van Diejen and 
V.~P.~Spiridonov, 
An elliptic Macdonald--Morris conjecture and multiple modular hypergeometric sums. Math. Res. Lett. {\bf 7} (2000), 729--746.
%

\bibitem{CG2007}
H.~Coskun and R.~A.~Gustafson, Well-poised Macdonald functions $W_\lambda$ and Jackson coefficients $\omega_\lambda$ on $BC_n$; 
in {\em Jack, Hall--Littlewood and Macdonald polynomials}, pp.127--155, Contemp. Math., {\bf 417}, Amer. Math. Soc., Providence, RI, 2006. 

\bibitem{FT1997}
I.~B.~Frenkel and V.~G.~Turaev, 
Elliptic solutions of the Yang--Baxter equation and modular hypergeometric functions; in
{\em The Arnold--Gelfand mathematical seminars}, pp.171--204, Birkh\"auser Boston, Boston, MA, 1997. 

\bibitem{GR2004} G.~Gasper and M.~Rahman, 
{\it Basic hypergeometric series}, Second edition, Encyclopedia of Mathematics and its Applications, 96. Cambridge University Press, Cambridge, 2004. xxvi+428 pp. 

\bibitem{Ito2008}
M.~Ito, 
A multiple generalization of Slater's transformation formula for a very-well-poised-balanced ${}_{2r}\psi_{2r}$ series. Q. J. Math. {\bf 59} (2008), 
221--235.

\bibitem{Ra2006}
E.~M.~Rains, $BC_n$-symmetric Abelian functions. Duke Math. J. {\bf 135} (2006), 99--180.

\bibitem{Ra2010}
E.~M.~Rains, Transformations of elliptic hypergeometric integrals, Ann. Math. {\bf 171} (2010), 169--243.

\bibitem{Ro2001}
H.~Rosengren,  
A proof of a multivariable elliptic summation formula conjectured by Warnaar. 
In: $q$-series with applications to combinatorics, number theory, and physics (Urbana, IL, 2000), 193--202, Contemp. Math., {\bf 291}, Amer. Math. Soc., Providence, RI, 2001.

\bibitem{Ru1997}
S.\,N.\,M.\,Ruijsenaars,
First order analytic difference equations and integrable quantum systems. J. Math. Phys. {\bf 38} (1997), 1069--1146.

\bibitem{Sc2007}
M.~Schlosser, Elliptic enumeration of nonintersecting lattice paths. J. Combin. Theory Ser. A {\bf 114} (2007), 505--521.

\bibitem{Sp2003}
V.~P.~Spiridonov, Theta hypergeometric integrals. Algebra i Analiz {\bf 15} (2003), 161--215; translation in St. Petersburg Math. J. {\bf 15} (2004), 929--967


\bibitem{SW2011}
V.~P.~Spiridonov and S.~O.~Warnaar, New multiple $_6\psi_6$ summation formulas and related conjectures. Ramanujan J. {\bf 25} (2011), 319--342.

\bibitem{W2002}
S.~O.~Warnaar, 
Summation and transformation formulas for elliptic hypergeometric series. Constr. Approx. {\bf 18} (2002), 479--502.

\end{thebibliography}

% Non-BibTeX users please use

\end{document}